\documentclass{amsart}
\usepackage{amsmath, amssymb, amscd, amsthm}

\newcommand{\CB}{{\rm CB}}
\newcommand{\CH}{{\rm CH}}
\newcommand{\CC}{{\rm CC}}

\newcommand{\QCH}{{\rm QCH}}

\newcommand{\lmod}[1]{{#1}\text{\rm --{\bf Mod}}}
\newcommand{\xra}[1]{\ensuremath{\xrightarrow{#1}}}
\newcommand{\xla}[1]{\ensuremath{\xleftarrow{#1}}}
\newcommand{\B}[1]{\ensuremath{\mathbb{#1}}}
\newcommand{\C}[1]{\ensuremath{\mathcal{#1}}}
\newcommand{\G}[1]{\ensuremath{\mathfrak{#1}}}
\newcommand{\eotimes}[1]{\ensuremath{\underset{#1}{\otimes}}}

\newcommand{\emod}[2]{\text{\bf mod}_{#1}\text{--}{#2}}
\newcommand{\eproj}[2]{\text{\bf proj}_{#1}\text{--}{#2}}
\newcommand{\efree}[2]{\text{\bf free}_{#1}\text{--}{#2}}

\newtheorem{thm}{Theorem}[section]
\newtheorem{cor}[thm]{Corollary}
\newtheorem{prop}[thm]{Proposition}
\newtheorem{lem}[thm]{Lemma}
\theoremstyle{definition}
\newtheorem{defn}[thm]{Definition}
\newtheorem{rem}[thm]{Remark}
\newtheorem{exm}[thm]{Example}

\numberwithin{equation}{section}

\title{Hopf--Hochschild (co)homology of module algebras}
\author{Atabey Kaygun}
\email{\tt akaygun@uwo.ca}
\address{Department of Mathematics\\ University of Western Ontario\\ 
London, Ontario N6A 5B7\\ Canada}

\begin{document}
\maketitle 


\section{Introduction}

Our goal in this paper is to define a version of Hochschild homology
and cohomology suitable for a class of algebras admitting compatible
actions of bialgebras, called ``module algebras''
(Definition~\ref{ModuleAlgebra}).  Our motivation lies in the
following problem: for an algebra $A$ which admits a module structure
over an arbitrary bialgebra $B$ compatible with its product structure,
the Hochschild or the cyclic bicomplexes associated with this algebra
need not be differential graded $B$--modules.  The obstruction which
prevents these complexes from being $B$--linear is trivial whenever
the bialgebra $B$ is cocommutative, as in the case of group rings and
universal enveloping algebras.  Yet the same obstruction is far from
being trivial if the underlying bialgebra is non-cocommutative.  In
the sequel, we will investigate how much of the Hochschild homology is
retained after dividing this obstruction out.  To this end, we will
construct a new differential graded $B$--module $\QCH_*(A,B,V)$
(Proposition~\ref{Quotient} and Definition~\ref{HopfHochschild}) for a
$B$--module algebra $A$ and a $B$--equivariant $A$--bimodule $V$
(Definition~\ref{EquivariantModule}).  We will define $HH^{\rm
Hopf}_*(A,B,V)$ the Hopf--Hochschild homology of $A$ with coefficients
in $V$ as the homology of the complex $k\eotimes{B}\QCH_*(A,B,V)$.  We
would like to point out that the same strategy worked remarkably well
in the case of cyclic cohomology of module coalgebras.  In
\cite{Kaygun:BialgebraCyclicK} we show that if we start with the
cocyclic bicomplex of a module coalgebra twisted by a stable
anti-Yetter--Drinfeld module, dividing the analogous obstruction
results in the Hopf cyclic complex of
\cite{Khalkhali:HopfCyclicHomology} which was an extension of the Hopf
cyclic cohomology of Connes and Moscovici
\cite{ConnesMoscovici:HopfCyclicCohomology}.

In the context of cyclic (co)homology and $K$--Theory, one of the most
commonly used tools dealing with module algebras has been ``crossed
product algebras'' (Definition~\ref{CrossedProductAlgebra}).  There is
a large body of work dealing with algebras admitting actions of
discrete groups and compact Lie groups, e.g.
\cite{Green:KTheoryOfCrossedProductAlgebras,
  Julg:KTheoryOfCrossedProductAlgebras,
  AtiyahSegal:EquivariantKTheory,
  Brylinski:CyclicHomologyAndEquivariantTheories,
  Kasparov:EquivariantKKTheory,
  GetzlerJones:CyclicHomologyOfCrossedProductAlgebras,
  BlockGetzlerJones:CyclicHomologyOfCrossedProductAlgebras} and
references therein, which utilizes this tool to its fullest extent.
Also, there have been successful attempts in defining equivariant
cyclic (co)homology and $K$--Theory for module algebras over Hopf
algebras \cite{BaajSkandalis:HopfEquivariantKKTheory,
  KhalkhaliAkbarpour:EquivariantCyclicHomologyII,
  NeshveyevTuset:HopfEquivariantKTheory} again by using crossed
product algebras.  Crossed product algebras enter in our picture in
Corollary~\ref{TorInterpretation} where we show that Hopf--Hochschild
homology can also be defined as a derived functor on the category of
representations of a crossed product algebra.

The last result we prove is the Morita invariance of the
Hopf--Hochschild homology and cohomology
(Theorem~\ref{EquivariantMorita} and Theorem~\ref{TwistedMorita}).
Our proof utilizes some additional tools from functor homology
\cite{McCarthy:CyclicHomologyOfExactCategories,
  Pirashvili:HochschildAndCyclicHomologyAsFunctorHomology}.  In doing
so, we observe that the category of representations of a crossed
product algebra is rather small for computing equivariant invariants.
However, the short-comings of this category can be overcome by using
``$B$--categories'' (Definition~\ref{EquivariantCategory}).  We refer
the reader to Remark~\ref{NaiveQuotient} for a more detailed analysis
on this subject.

Here is the plan of this paper: In
Section~\ref{HopfHochschildHomology} we give the basic definition of
Hopf--Hochschild complex of module algebras with coefficients in
equivariant bimodules.  We also point out the connections between
Hopf--Hochschild homology and Hopf cyclic cohomology
(Remark~\ref{Connections}).  In Section~\ref{HopfHochschildCohomology}
we define Hopf--Hochschild cohomology of a $B$--module algebra and
calculate it for lower dimensions.  We also give a derived functor
interpretation of the Hopf--Hochschild cohomology in terms of crossed
product algebras.  In Section~\ref{HopfHochschildHomologyRevisited} we
extend the derived functor interpretation to Hopf--Hochschild
homology.  Section~\ref{CategoricalAlgebra} and
Section~\ref{EquivariantCategories} contain technical results needed
toward proving Morita invariance of Hopf--Hochschild (co)homology in
Section~\ref{Morita} in its full generality.  In
Section~\ref{TwistedBifunctors} we develop a generalized twisting
method for coefficient bifunctors or bimodules by using
Yetter-Drinfeld modules similar to the method of twisting developed in
\cite{Khalkhali:HopfCyclicHomology}.  In this last Section, we also
prove Morita invariance for twisted Hopf--Hochschild (co)homology.

Throughout this paper, we assume $k$ is an arbitrary field and $B$ is
an associative/coassociative unital/counital bialgebra, or a Hopf
algebra with an invertible antipode whenever it is necessary.  All
tensor products are taken over $k$ unless it is stated otherwise
explicitly.

We would like to thank S. Witherspoon for her helpful comments and
corrections.

\section{Hopf--Hochschild homology}
\label{HopfHochschildHomology}

\begin{defn}\label{ModuleAlgebra}
An algebra $A$ is called a left $B$--module algebra if $A$ is a
$B$--module and
\begin{align*}
b(a_1a_2) = b_{(1)}(a_1)b_{(2)}(a_2)
\end{align*}
for any $b\in B$ and $a_1,a_2\in A$.  If $A$ is unital, we also assume
$b(1_A)=\varepsilon(b)1_A$ where $\varepsilon$ is the counit of $B$.
\end{defn}

\begin{defn}\label{EquivariantModule}
Let $A$ be a $B$--module algebra.  An $A$--module $V$ is called a
$B$--equivariant $A$--module if $V$ is both an $A$--module and
$B$--module and one also has
\begin{align*}
b(av) = b_{(1)}(a)(b_{(2)}(v))
\end{align*}
for any $a\in A$ and $b\in B$.
\end{defn}

\begin{exm}
Let $B=k[G]$ be the group algebra of a discrete group $G$.  Then an
algebra $A$ is a $k[G]$--module algebra iff it is a $G$--algebra.
\end{exm}

\begin{exm}
Let $B=U(\G{g})$ is the universal enveloping algebra of Lie algebra
$\G{g}$. Then an algebra $A$ is a $U(\G{g})$--module algebra iff it
admits an action of $\G{g}$ by derivations.
\end{exm}

\begin{exm}
Let $B$ be an arbitrary Hopf algebra and $A$ be an algebra which is
also a $B$--bimodule.  Then there is a natural action of $B$ on $A$
called the adjoint action which makes $A$ into a $B$--module algebra
and any $A$--module which is also a $B$--module $V$ into a
$B$--equivariant $A$--module.  The adjoint action is defined as
\begin{align*}
ad_b(a) = b_{(1)}a S(b_{(2)})
\end{align*}
for any $b\in B$ and $a\in A$.  One can easily see that
\begin{align*}
ad_b(a_1a_2) = b_{(1)}a_1a_2 S(b_{(2)})
 = b_{(1)}a_1 S(b_{(2)}) b_{(3)}a_2S(b_{(4)})
 = ad_{b_{(1)}}(a_1) ad_{b_{(2)}}(a_2)
\end{align*}
for any $b\in B$ and $a_1,a_2\in A$.  Similarly
\begin{align*}
b(av) = b_{(1)}a S(b_{(2)}) b_{(3)} v
  = ad_{b_{(1)}}(a)(b_{(2)}v)
\end{align*}
for any $a\in A$, $b\in B$ and $v\in V$.
\end{exm}

\begin{defn}
Given an algebra $A$ and a $A$--bimodule $V$, we will use the notation
$\CH_*(A,V)$ to denote graded module $\bigoplus_{n\geq 0} A^{\otimes
n} \otimes V$ with structure morphisms
\begin{align*}
\partial_j(a_1\otimes\cdots\otimes a_n\otimes v)
 = \begin{cases}
   (\cdots\otimes a_{j+1}a_{j+2}\otimes\cdots\otimes v)
         & \text{ if } 0\leq j< n-1\\
   (a_1\otimes\cdots\otimes a_nv)
         & \text{ if } j=n-1\\
   (a_2\otimes\cdots\otimes a_n\otimes va_1)
         & \text{ if } j=n
   \end{cases}
\end{align*}
which makes $\CH_*(A,V)$ into a pre-simplicial module.  The
differential graded module with the differentials
\begin{align*}
d^{\CH}_n = \sum_{j=0}^n (-1)^j\partial_j
\end{align*}
corresponding to this pre-simplicial module is also denoted by
$\CH_*(A,V)$, and is called the Hochschild complex of $A$ with
coefficients in the $A$--bimodule $V$.
\end{defn}

From this point on, we will assume $A$ is a $B$--module algebra and $V$
is a $B$--equivariant $A$--bimodule unless it is stated otherwise
explicitly.

\begin{rem}
$B$ as an algebra acts on $\CH_*(A,V)$ diagonally as
\begin{align*}
L_b(a_1\otimes\cdots\otimes a_n\otimes v)
 = b_{(1)}(a_1)\otimes\cdots\otimes b_{(n)}(a_n)\otimes b_{(n+1)}(v)
\end{align*}
which makes $\CH_*(A,V)$ into a graded $B$--module but NOT a
differential graded $B$--module since
\begin{align*}
\partial_n L_b(a_1\otimes\cdots\otimes a_n\otimes v)
 = &   b_{(2)}(a_2)\otimes\cdots\otimes b_{(n)}(a_n)\otimes
       b_{(n+1)}(v) b_{(1)}(a_1)\\
\neq & b_{(1)}(a_2)\otimes\cdots\otimes b_{(n-1)}(a_n)\otimes
       b_{(n)}(v)b_{(n+1)}(a_1)\\
 = & L_b \partial_n(a_1\otimes\cdots\otimes a_n\otimes v)
\end{align*}
unless $b\in ker(\delta)$ where $\delta = (id_2-\tau_2)\Delta$.  This
means, although the $B$--structure on $A$ does extend to a graded
$B$--module structure on the ordinary Hochschild complex $\CH_*(A,V)$,
it extends to a differential graded $B$--module structure on
$\CH_*(A,V)$ {\bf if} $B$ is cocommutative, for instance when $B$ is a
group ring or a universal enveloping algebra.  The obstruction which
prevents $\CH_*(A,V)$ from being a differential graded $B$--module is
the subcomplex generated by images of the the commutators
$[L_x,d^\CH_*]$ where $L_x$ is the $k$--linear endomorphism of
$\CH_*(A,V)$ coming from the diagonal action of $x\in B$ on
$\CH_*(A,V)$.  Now one can ask the following question: what happens if
we force these differential graded $k$--modules to become differential
graded $B$--modules by dividing out this obstruction?  This is what we
are going to do with Definition~\ref{TheIdeal} and
Proposition~\ref{Quotient} for the ordinary Hochschild complex.  In
the sequel, we investigate homological consequences of this operation.
\end{rem}

Let $(\C{C}_*,d^\C{C}_*)$ be a differential graded $k$--module and let
$n\in\B{N}$.  Then we define $\C{C}_*[+n]$ the $n$-fold suspension of
$\C{C}_*$ as the differential graded graded $k$--module
$\C{C}_m[+n]=\C{C}_{m+n}$ with differentials
$d^\C{C}_m[+n]=d^\C{C}_{m+n}$ for any $m\in\B{Z}$.  One can similarly
define $\C{C}_*[-n]$ for any $n\in\B{N}$.  Note that
$H_{m\pm n}(\C{C}_*[\pm n])= H_m(\C{C}_*)$.

\begin{lem}\label{DGM}
For any $b\in B$, there is a morphism of differential graded
$k$--modules of the form
\begin{align*}
  \CH_*(A,V)[+1]\xra{[L_b,\partial_{*+1}]}\CH_*(A,V)
\end{align*}
Moreover, $[L_b,\partial_{*+1}]$ is null-homotopic for any $b\in B$.
\end{lem}

\begin{proof}
For any $b\in B$ we consider
\begin{align*}
\partial_j[L_b,\partial_{n+1}] & ({\bf a}\otimes v)
 = - [L_b,\partial_j]\partial_{n+1}({\bf a}\otimes v)
    + [L_b,\partial_j\partial_{n+1}]({\bf a}\otimes v)\\
 =& \begin{cases}
     [L_b,\partial_n\partial_j]({\bf a}\otimes v) 
              & \text{ if } 0\leq j\leq n-1\\
     - [L_b,\partial_n]\partial_{n+1}({\bf a}\otimes v)
     + [L_b,\partial_n\partial_n]({\bf a}\otimes v) 
              & \text{ if } j=n
    \end{cases}\\
 =& \begin{cases}
     [L_b,\partial_n]\partial_j({\bf a}\otimes v) 
              & \text{ if } 0\leq j\leq n-1\\
     - [L_b,\partial_n]\partial_{n+1}({\bf a}\otimes v)
     + [L_b,\partial_n]\partial_n({\bf a}\otimes v) 
              & \text{ if } j=n
    \end{cases}
\end{align*}
by using the fact that $[L_b,\partial_j]({\bf a}\otimes v)=0$ for any
$0\leq n$, $0\leq j\leq n-1$ and for any $({\bf a}\otimes v)$ from
$\CH_n(A,V)$.  This result immediately implies 
\begin{align*}
d^{\CH}_n[L_b,\partial_{n+1}]=[L_b,\partial_n]d^{\CH}_n
\end{align*}
The null-homotopy $\CH_n(A,V)\xra{s_n}\CH_{n+1}(A,V)[+1]$ is given
by $s_n=(-1)^{n-1}L_b$ for any $n\geq 0$.
\end{proof}

\begin{defn}\label{TheIdeal}
We define a graded $B$--submodule of $\CH_*(A,V)$ as
\begin{align*}
J_*(A,B,V)
 = \sum_{b\in B} im([L_b,\partial_{*+1}])
\end{align*}
\end{defn}

\begin{prop}\label{Quotient}
 We define a new graded $B$--module $\QCH_*(A,B,V)$ as the quotient
 graded $B$--module 
 \begin{align*}
   \QCH_*(A,B,V):=\CH_*(A,V)/J_*(A,B,V)
 \end{align*}
 Then $\QCH_*(A,B,V)$ is also a differential graded $B$--module.
\end{prop}

\begin{proof}
Since each $im([L_b,\partial_{*+1}])$ is a differential graded
$k$--submodule of $\CH_*(A,V)$, the submodule $J_*(A,B,V)$ is also a
differential graded $k$--submodule of $\CH_*(A,V)$.  This fact implies
$\QCH_*(A,B,V)$ is too a differential graded $k$--module.  Moreover,
$J_*(A,B,V)$ is a graded $B$--submodule of $\CH_*(A,V)$ since
\begin{align*}
L_x[L_b,\partial_{n+1}](v\otimes{\bf a})
 = -[L_x,\partial_{n+1}]L_b(v\otimes{\bf a})
   +[L_{xb},\partial_{n+1}](v\otimes{\bf a})
\end{align*}
for any $x,b\in B$, $n\geq 0$ and $(v\otimes{\bf a})$ from
$\CH_n(A,V)[+1]$.  In order $\QCH_*(A,B,V)$ be a differential graded
$B$--module, we must show that $[L_b,d^\CH_*]\equiv 0$ on
$\QCH_*(A,B,V)$ for any $b\in B$.  This is equivalent to saying that
\begin{align*}
[L_b,d^\CH_n](v\otimes{\bf a})
  = (-1)^n[L_b,\partial_n](v\otimes{\bf a})
\end{align*}
must be in $J_{n-1}(B,B,V)$ for any $(v\otimes{\bf a})$ from
$\CH_n(A,V)$ and for any $b\in B$ is in $J_{n-1}(A,B,V)$ which is true
by definition.
\end{proof}

\begin{defn}\label{HopfHochschild}
Assume $B$ is a bialgebra.  For a $B$--module algebra $A$ and a
$B$--equivariant $A$--module $V$, we define Hopf--Hochschild
homology of $A$ with coefficients in $V$ as the homology of the
differential graded $k$--module ${}_B\QCH_*(A,B,V)$ which is defined
as $k\eotimes{B}\QCH_*(A,B,V)$.  In other words,
\begin{align*}
HH^{\rm Hopf}_n(A,V) := H_n({}_B\QCH_*(A,B,V))
\end{align*}
for any $n\geq 0$.
\end{defn}

\begin{rem}
  Observe that if $B$ is cocommutative, the differential graded
  $k$--module $\QCH_*(A,B,V)$ is equal to the ordinary Hochschild
  complex $\CH_*(A,V)$.  In an (un)likely case when $B$ is both
  cocommutative and semi-simple (such as $B=k[G]$ where $G$ is a
  finite group and $char(k)$ does not divide $|G|$) then one has an
  isomorphism of the form $HH^{\rm Hopf}_*(A,V)\cong
  k\eotimes{B}HH_*(A,V)$.  For example, if $B=k$ the Hopf--Hochschild
  homology is the same as the ordinary Hochschild homology.
\end{rem}

\begin{rem}\label{Connections}
  Assume $A$ is an associative, but not necessarily unital
  $k$--algebra.  Apart from the ordinary Hochschild complex of $A$,
  there are several other different differential graded $k$--modules
  one can associate with $A$:
  \begin{enumerate}




  \item Connes' complex $\CC^\lambda_*(A)$ which is defined as the
    cyclic coinvariants of the ordinary Hochschild complex.

  \item The positive, negative and periodic cyclic bicomplexes
    $\CC_*(A)$, ${\rm CN}_*(A)$ and $\CC^{\rm
    per}_*(A)$.
 
  \item The mixed complex ${\rm CM}_*(A)$, which is also
    referred as ``the $(b,B)$--complex'' which also has two other
    variations: the negative mixed complex ${\rm CM}_*^-(A)$
    and the periodic mixed complex ${\rm CM}_*^{\rm per}(A)$.

  \end{enumerate}
  As before, the cyclic bicomplexes and the mixed complexes are graded
  $B$--modules but are not necessarily differential graded
  $B$--modules.  The obstruction to extending the graded $B$--module
  structure to a differential graded $B$--module structure stems from
  the fact that the cyclic permutations and the diagonal $B$--action
  on the tensor powers of $A$ do not necessarily commute.  We will
  investigate the consequences of the operation of dividing out this
  obstruction on the cyclic complexes we mentioned above in a
  different paper \cite{Kaygun:BivariantHopf} in a more general set-up
  where the complexes are twisted by some coefficient module.  We
  would like to point out that the obstruction which prevents the
  cyclic bicomplex from being a differential graded $B$--module is a
  larger differential graded submodule in the sense that the
  ``Hochschild subcomplex'' or the ``$b$--subcomplex'' of the Hopf
  cyclic bicomplex is a quotient of the Hopf--Hochschild complex we
  define here.
\end{rem}

\section{Hopf--Hochschild cohomology}
\label{HopfHochschildCohomology}

\begin{defn}\label{CrossedProductAlgebra}
We define $A^e\rtimes B$ as $A^e\otimes B=A\otimes A^{op}\otimes B$
with the multiplication
\begin{align*}
(a_1\otimes a_1'\otimes b^1)(a_2\otimes a_2'\otimes b^2)
 = (a b^1_{(1)}(a_2)\otimes b^1_{(3)}(a_2')a_1'\otimes b^1_{(2)}b^2)
\end{align*}
for any $(a_1\otimes a_1'\otimes b^1)$ and $(a_2\otimes a_2'\otimes
b^2)$ from $A^e\rtimes B$.
\end{defn}

\begin{lem}
$A^e\rtimes B$ is a unital associative algebra.
\end{lem}

\begin{proof}
For associativity, one must consider
\begin{align*}
((a_1\otimes a_1'\otimes b^1) & (a_2\otimes a_2'\otimes b^2))
(a_3\otimes a_3'\otimes b^3)\\
 & = (a_1 b^1_{(1)}(a_2)\otimes b^1_{(3)}(a_2')a_1'\otimes b^1_{(2)}b^2)
     (a_3\otimes a_3'\otimes b^3)\\
 & = a_1 b^1_{(1)}(a_2b^2_{(1)}(a_3))\otimes 
     b^1_{(3)}(b^2_{(3)}(a_3')a_2')a_1'\otimes 
     b^1_{(2)}b^2_{(2)}b^3\\
 & = (a_1\otimes a_1'\otimes b^1) 
     ((a_2\otimes a_2'\otimes b^2)(a_3\otimes a_3'\otimes b^3))
\end{align*}
for any $(a_i\otimes a_i'\otimes b_i)$ from $A^e\rtimes B$ for
$i=1,2,3$. 
\end{proof}

\begin{lem}
$V$ is a $B$--equivariant $A$--bimodule iff $V$ is a $A^e\rtimes
B$--module.
\end{lem}

\begin{proof}
$V$ is a $B$--equivariant $A$--bimodule iff one has
\begin{align*}
(a_1\otimes a_2\otimes b)((a_1'\otimes a_2'\otimes b')(v))
 & = a_1 b(a_1' b'(v) a_2') a_2
   = a_1 b_{(1)}(a_1') (b_{(2)} b')(v) b_{(3)}(a_2') a_2\\
 & = (a_1b_{(1)}(a_1')\otimes b_{(3)}(a_2')a_2)\cdot(b_{(2)}b')(v)\\
 & = ((a_1\otimes a_2\otimes b)(a_1'\otimes a_2'\otimes b'))
     \cdot v
\end{align*}
for any $(a_1\otimes a_2\otimes b)$ and $(a_1'\otimes a_2'\otimes b')$
from $A^e\rtimes B$ and $v\in V$, i.e. $V$ is a $A^e\rtimes B$--module.
\end{proof}

\begin{defn}
Define a differential graded $k$--module $\CB_*(A)$ by letting
$\CB_n(A)=A^{\otimes n+2}$ for any $n\geq 0$.  Then we define
pre-simplicial structure morphisms
\begin{align*}
\partial_j(a_0\otimes\cdots\otimes a_{n+1})
 = (a_0\otimes\cdots\otimes a_j a_{j+1}\otimes\cdots\otimes a_{n+1})
\end{align*}
for any $n\geq 1$ and $0\leq j\leq n$ and then define the
differentials as $d^\CB_n=\sum_{j=0}^n(-1)^j\partial_j$ for any $n\geq
1$.
\end{defn}

\begin{lem}
$\CB_*(A)$ is a differential graded $A^e\rtimes B$--module.
\end{lem}

\begin{proof}
$\CB_*(A)$ is a differential graded $A$--bimodule with the $A$--action
defined as
\begin{align*}
a(a_0\otimes\cdots\otimes a_{n+1})
 & = (aa_0\otimes\cdots\otimes a_{n+1})  & 
(a_0\otimes\cdots\otimes a_{n+1})a 
 & = (a_0\otimes\cdots\otimes a_{n+1}a)
\end{align*}
for any $a\in A$ and $(a_0\otimes\cdots\otimes a_{n+1})$ from
$\CB_n(A)$.  There is also a left $B$--module structure defined for
$b\in B$ and $(a_0\otimes\cdots\otimes a_{n+1})$ from $\CB_n(A)$ as
\begin{align*}
b(a_0\otimes\cdots\otimes a_{n+1})
 & = (b_{(1)}a_0\otimes\cdots\otimes b_{(n+2)}a_{n+1})
\end{align*}
This makes $\CB_*(A)$ into a left $A^e\rtimes B$--module since one has
\begin{align*}
b(a\otimes a') & (a_0\otimes\cdots\otimes a_{n+1})\\
 & = (b_{(1)}(aa_0)\otimes b_{(2)}a_1\otimes\cdots\otimes 
        \otimes b_{(n+1)}a_n\otimes b_{(n+2)}(a_{n+1}a'))\\
 & = (b_{(1)}(a)b_{(2)}(a_0)\otimes b_{(3)}a_1\otimes\cdots\otimes 
      b_{(n+3)}(a_{n+1}) b_{(n+4)}(a'))\\
 & = (b_{(1)}(a)\otimes b_{(3)}(a'))b_{(2)}(a_0\otimes\cdots\otimes a_{n+1})
\end{align*}
$b\in B$, $(a\otimes a')$ from $A^e$ and $(a_0\otimes\cdots\otimes
a_{n+1})$ from $\CB_n(A)$.
\end{proof}

\begin{defn}
Let $A$ be a $B$--module algebra and let $V$ be an $B$--equivariant
$A$--bimodule.  Define Hopf--Hochschild cochain complex of $A$
with coefficients in $V$ as 
\begin{align*}
  \CH_{\rm Hopf}^*(A,V) := {\rm Hom}_{A^e\rtimes B}(\CB_*(A),V)
\end{align*}
$HH_{\rm Hopf}^*(A,V)$ Hopf--Hochschild cohomology of $A$ with
coefficients in $V$ is defined to be the cohomology of this cochain
complex.
\end{defn}

\begin{thm}\label{ExtInterpretation}
  Let $A$ be an $H$--unital projective $B$--module algebra and assume
  $V$ is an arbitrary $B$--equivariant $A$--bimodule.  Then one has
  isomorphisms of the form
  \[ HH_{\rm Hopf}^{n+1}(A,V) \cong {\rm Ext}_{A^e\rtimes B}^n(\Omega(A),V) \]
  for any $n\geq 1$.
\end{thm}

\begin{proof}
  We have a short exact sequence of $A^e\rtimes B$-modules of the form
  \[ 0 \xra{} \Omega(A) \xra{} A^e \xra{\mu_A} A \xra{} 0 \] where
  $\mu_A$ denotes the multiplication morphism.  Since $A$ is
  $B$--projective, the brutal truncation $\CB_{*>0}(A)$ of $CB_*(A)$
  at $n=0$ is a $A^e\rtimes B$--projective resolution of $\Omega(A)$
  the kernel of the multiplication map.  The result follows.
\end{proof}

\begin{defn}
Assume $A$ is a $B$--module algebra and $V$ is a $B$--equivariant
$A$--module.  A morphism of $k$--modules $A\xra{D}V$ is called a
$V$--valued derivation on $A$ iff
\[ D(aa') = D(a) a' + a D(a') \]
for any $a,a'\in A$.  The same derivation is called a $V$--valued
$B$--equivariant derivation if
\[ D(ba) = b D(a) \] 
for any $a\in A$ and $b\in B$.  The $k$--module of $V$--valued and
$V$--valued $B$--equivariant derivations on $A$ are denoted by ${\rm
Der}(A,V)$ and ${\rm Der}_B(A,V)$ respectively.
\end{defn}

Let $v\in V$ be fixed and consider the $k$--module morphism
$A\xra{[v,\cdot]}V$ defined by $[v,a]=va-av$ for any $a\in A$.  Then
\begin{align*}
 [v,aa'] = & vaa' - aa'v = vaa' - ava' + ava' - aa'v
         = [v,a]a' + a[v,a']
\end{align*}
for any $a,a'\in A$ and $b\in B$.  This means elements of the form
$[v,\cdot]$ from ${\rm Hom}_k(A,V)$ are $V$--valued derivations on
$A$.  However if we were to require $[v,\cdot]$ to be a
$B$--equivariant derivation, then we need to have
\begin{align*}
 [v,b(a)] = & vb(a) - b(a)v 
          = b_{(1)}(v) b_{(2)}(a) - b_{(1)}(a)b_{(2)}(v)
          = b[v,a]
\end{align*}
for any $a\in A$, $b\in B$.  In case $V$ is a trivial $B$--module via
the counit, i.e. $b(v) = \varepsilon(b)v$ then the condition above is
satisfied. 

\begin{prop}
  Let $A$ be an $B$-module algebra and $V$ be an $B$--equivariant
  $A$--module as before.  Then one has 
  \begin{align*}
    HH^0_{\rm Hopf}(A,V)\cong  ({}^B V)^{{\rm Lie}(A)}
    \text{\ \  and\ \ } 
    HH^1_{\rm Hopf}(A,V) \cong {\rm Der}_B(A,V)/[A,{}^B V]
\end{align*}
\end{prop}

\begin{proof}
If we consider $f\in{\rm Hom}_{A^e\rtimes B}(\CB_n(A),V)$ we see that
\begin{align*}
f(a_0\otimes\cdots\otimes a_{n+1})
 = a_0f(1\otimes a_1\otimes\cdots\otimes a_n\otimes 1)a_{n+1}
\end{align*}
Since $b(1_A)=\varepsilon(b)$ for all $b\in B$, we see that $\CH_{\rm
Hopf}^0(A,V)\cong {\rm Hom}_B(k,V)\cong {}^B V$ where ${}^B V$ is
defined as the submodule $\{v\in V|\ b(v)=\varepsilon(b)v,\ b\in B\}$
of $V$.  The proof that $\CH^1_*(A,V)\cong {\rm Hom}_B(A,V)$ is
similar.  Then $v\in HH^0_{\rm Hopf}(A,V)$ iff
\begin{align}\label{EquivHoch0}
d_\CH^0(v)(1\otimes a\otimes 1) = v(a\otimes 1) - v(1\otimes a)
 = a v - v a = 0
\end{align}
i.e. $[a,v]=0$ for any $a\in A$ which is the same as the invariants of
the adjoint action of $\text{Lie}(A)$ on ${}^B V$,
i.e. $({}^B V)^{\text{Lie}(A)}$.  Similarly, $f\in ker(d_{\CH}^1)$ iff
\begin{align*}
(d_{\CH}^1 f)(1\otimes a\otimes a'\otimes 1)
 = & f(a\otimes a'\otimes 1) - f(1\otimes aa'\otimes 1) 
  + f(1\otimes a\otimes a')\\
 = & a f(1\otimes a'\otimes 1) - f(1\otimes aa'\otimes 1) 
  + f(1\otimes a\otimes 1) a'
 = 0
\end{align*}
for any $(1\otimes a\otimes a'\otimes 1)$ from $\CB_2(A)$. In other
words $f\in ker(d_{\CH}^1)$ iff $D_f(a) := f(1\otimes a\otimes 1)$ is
a derivation.  Moreover, since $f$ is $B$--equivariant, so is $D_f$.
Equation \ref{EquivHoch0} tells us that the image of $d_{\CH}^0$
consists of the elements of the form $[a,v]$ where $v\in {}^B V$ and
$a\in A$.  In other words $HH_{\rm Hopf}^1(A,V)\cong {\rm
  Der}_B(A,V)/[A,{}^B V]$.
\end{proof}

\section{Hopf--Hochschild homology revisited}
\label{HopfHochschildHomologyRevisited}

\begin{defn}
Assume $B$ is a Hopf algebra and let $U$ be a right $B$--module.  Then
one can think of $U$ as a left module via the action $b\cdot u :=
u\cdot S(b)$ for any $u\in U$ and $b\in B$.  We denote the new module
by $U^{op}$.
\end{defn}

\begin{thm}\label{MainIso}
Assume $B$ is a Hopf algebra with an invertible antipode.  Let $A$ be
a $B$--module algebra and $V$ be a $B$--equivariant $A$--bimodule.
Then ${}_B \QCH_*(A,B,V)$ is isomorphic to $V^{op} \eotimes{A^e\rtimes
B}\CB_*(A)$ as differential graded $k$--modules.
\end{thm}

\begin{proof}
Define a morphism of graded modules 
\begin{align*}
\CH_*(A,V)\xra{\varphi_*}V^{op}\eotimes{A^e\rtimes B}\CB_*(A)
\end{align*}
by letting
\begin{align*}
\varphi_n(a_1\otimes\cdots\otimes a_n\otimes v)
 = v\eotimes{A^e\rtimes B}(1\otimes a_1\otimes\cdots\otimes a_n\otimes 1)
\end{align*}
for any $(a_1\otimes\cdots\otimes a_n\otimes v)$ from $\CH_n(A,V)$.
Notice that 
\begin{align*}
\varphi L_b(a_1\otimes\cdots\otimes a_n\otimes v)
& = v S^{-1}(b_{(2)})\eotimes{A^e\rtimes B} 
    (1\otimes L_{b_{(1)}}(a_1\otimes\cdots\otimes a_n)\otimes 1)\\
& = v S^{-1}(b_{(2)})\eotimes{A^e\rtimes B} 
    L_{b_{(1)}}(1\otimes a_1\otimes\cdots\otimes a_n\otimes 1)\\
& = \varepsilon(b)v\eotimes{A^e\rtimes B} 
    (1\otimes a_1\otimes\cdots\otimes a_n\otimes 1)
\end{align*}
by observing $x(1_A)=\varepsilon(x)$ for any $x\in B$.  Therefore $\varphi_*$
factors as 
\begin{align*}
\CH_*(A,V)\xra{q'_*}{}_B\CH_*(A,V)\xra{\varphi'_*}
    V^{op}\eotimes{A^e\rtimes B}\CB_*(A)
\end{align*}
where $q'_*$ is the canonical quotient taking $\CH_*(A,V)$ to the graded
module of $B$--coinvariants ${}_B\CH_*(A,V)$.  However, for any $b\in B$
we also have
\begin{align*}
\varphi_n[L_b,\partial_{n+1}]
 = \varphi_n L_b\partial_{n+1} - \varphi_n \partial_{n+1}L_b
 = \varphi_n \partial_{n+1} (\varepsilon(b)-L_b)
\end{align*}
Together with the fact that $\varphi_n \partial_{n+1} =
(id_V\otimes\partial_0)\varphi_{n+1}$
we see 
\begin{align*}
\varphi_n[L_b,\partial_{n+1}]
 = (id_V\otimes\partial_0)\varphi_{n+1}(\varepsilon(b)-L_b) = 0
\end{align*}
for any $b\in B$.  This implies we get an extension
\begin{align*}
\begin{CD}
\CH_*(A,V)  @>{q'_*}>> {}_B\CH_*(A,V)
              @>{\varphi'_*}>> V^{op}\eotimes{A^e\rtimes B}\CB_*(A)\\
@V{q_*}VV     @V{{}_B q_*}V{\cong}V           @|\\
\QCH_*(A,B,V) @>>{q'_*}> {}_B\QCH_*(A,B,V)
              @>>{\varphi''_*}> V^{op}\eotimes{A^e\rtimes B}\CB_*(A)\\
\end{CD}
\end{align*}
Since $\varphi_*$ is a morphism of differential graded $k$--modules so
is $\varphi''_*$ because both $q_*$ and $q_*'$ are morphisms of
differential graded $k$--modules too.  Notice that ${}_Bq_*$ is an
isomorphism of graded $k$--modules since $q'_*$ factors through
$\QCH_*(A,B,V)$.  Also, observe that $V^{op}\eotimes{A^e\rtimes
B}\CB_n(A)$ is generated by elements of the form
$(v\eotimes{A^e\rtimes B}(1\otimes a_1\otimes\cdots\otimes a_n\otimes
1))$ since
\begin{align*}
(v\eotimes{A^e\rtimes B} & (a_0\otimes a_1\otimes\cdots\otimes 
  a_n\otimes a_{n+1}))\\
&=(v\eotimes{A^e\rtimes B}(a_0\otimes a_{n+1})(1\otimes a_1\otimes
    \cdots\otimes a_n\otimes 1))\\
&=(a_{n+1}va_0\eotimes{A^e\rtimes B}
     (1\otimes a_1\otimes \cdots\otimes a_n\otimes 1))
\end{align*}
for some $v\in V$ and $a_i\in A$.  This implies $\varphi_*''$ is an
epimorphism.  Now define a section
\[ V^{op}\eotimes{A^e\rtimes B}\CB_*(A)\xra{s_*}{}_B\QCH_*(A,B,V) \] 
by letting
\begin{align*}
s_n(v\eotimes{A^e\rtimes B}(a\otimes a_1\otimes\cdots\otimes a_n\otimes a'))
=& [a_1\otimes\cdots\otimes a_n\otimes a'v a]
\end{align*}
for any $(v\eotimes{A^e\rtimes B}(1\otimes a_1\otimes\cdots\otimes
a_n\otimes 1))$ from $V^{op}\eotimes{A^e\rtimes B}\CB_*(A)$.  The
section $s_*$ is well-defined since
\begin{align*}
s_n & \left(v(a\otimes a'\otimes b)\eotimes{A^e\rtimes B}
            (1\otimes a_1\otimes\cdots\otimes a_n\otimes 1)\right)\\
=&s_n\left(S(b)(a'va)\eotimes{A^e\rtimes B}
    (1\otimes a_1\otimes\cdots\otimes a_n\otimes 1)\right)\\
=& [a_1\otimes\cdots\otimes a_n\otimes S(b)(a'va)]
=  [S(b_{(2)})b_{(3)}(a_1\otimes\cdots\otimes a_n)\otimes
    S(b_{(1)})(a'va)]\\
=& [b(a_1\otimes\cdots\otimes a_n)\otimes a' v a]
= s_n\left(v\eotimes{A^e\rtimes B}
            (a\otimes b(a_1\otimes\cdots\otimes a_n)\otimes a')\right)\\
=& s_n\left(v\eotimes{A^e\rtimes B}(a\otimes a'\otimes b)
            (1\otimes a_1\otimes\cdots\otimes a_n\otimes 1)\right)
\end{align*}
for any $(v\eotimes{A^e\rtimes B}(1\otimes a_1\otimes\cdots\otimes
a_n\otimes 1))$ from $V^{op}\eotimes{A^e\rtimes B}\CB_*(A)$ and
$(a\otimes a'\otimes b)$ from $A^e\rtimes B$.  We used the notation
$[\Psi]$ to denote the class of an element $\Psi\in \QCH_*(A,V)$ in
${}_B\QCH_*(A,V)$.  One can easily see that $s_*\varphi'_*=id_*$ which
means $\varphi''_*$ is also a monomorphism.  The result follows.
\end{proof}

\begin{cor}\label{Tor0}
Let $B$ be a Hopf algebra with an invertible antipode.  Assume $A$ is a
$B$--module algebra and $V$ is a left $B$--equivariant
$A$--bimodule. Then
\begin{align*}
HH^{\rm Hopf}_0(A,V^{op})\cong {}_BV/[A,{}_BV]
\end{align*}
where ${}_BV:= k\eotimes{B}V \cong A^e\eotimes{A^e\rtimes B}V$.
\end{cor}

\begin{cor}\label{TorInterpretation}
Let $B$ be a Hopf algebra with an invertible antipode.  Assume $A$
is an $H$--unital projective $B$--module algebra and $V$ is a
$B$--equivariant $A$--bimodule.  Then one has isomorphisms of the form
\[ HH^{\rm Hopf}_{n+1}(A,V) \cong {\rm Tor}^{A^e\rtimes B}_n(\Omega(A),V^{op}) \]
for any $n\geq 1$.
\end{cor}

\begin{proof}
  Consider the short exact sequence of $A^e\rtimes B$-modules
  \[ 0\xra{}\Omega(A)\xra{}A^e\xra{\mu_A}A\xra{}0 \] where $\mu_A$ is
  the multiplication on $A$.  If $A$ is $H$--unital and $A$ is
  $B$--projective then the brutal truncation $\CB_{*>0}(A)$ of
  $CB_*(A)$ at $n=0$ is a $A^e\rtimes B$--projective resolution of
  $\Omega(A)$.  The rest of the proof is similar to that of
  Theorem~\ref{ExtInterpretation}.
\end{proof}

\begin{rem}
After Theorem~\ref{MainIso} and Corollary~\ref{TorInterpretation},
one can see that there is another possible definition of the
Hopf--Hochschild complex of a $B$--module algebra $A$.  Namely, one
can use the differential graded $k$--module
\begin{align*}
V\eotimes{A^e\rtimes B}\CB_*(A)
\end{align*}
to define the Hopf--Hochschild homology of $A$ for a right
$B$--equivariant $A$--bimodule $V$.
\end{rem}

\section{Categorical algebra and cofinality}
\label{CategoricalAlgebra}

\begin{defn}
  A small category $\C{C}$ is called $k$--linear if for each $X,Y\in
  Ob(\C{C})$, the Hom object ${\rm Hom}_\C{C}(X,Y)$ is a $k$--module
  and the composition maps
  \[ {\rm Hom}_\C{C}(Y,Z)\times {\rm Hom}_\C{C}(X,Y)\xra{} {\rm Hom}_\C{C}(X,Z) \]
  are $k$--bilinear for any $X,Y,Z\in Ob(\C{C})$.  A functor
  $\C{C}\xra{F}\C{C}'$ between two $k$--linear categories is called a
  $k$--linear functor if the map
  \[ {\rm Hom}_\C{C}(X,Y) \xra{F_{X,Y}} {\rm Hom}_{\C{C}'}(X,Y) \]
  is a morphism of $k$--modules for any $X,Y$ in $Ob(\C{C})$.
\end{defn}

\begin{defn}
  A $k$--linear bifunctor on $\C{C}$ with values in another
  $k$--linear category $\C{C}'$ is just a functor of the form
  $\C{C}\times\C{C}^{op}\xra{\C{H}}\C{C}'$.
\end{defn}

\begin{defn}
  Let $\C{C}$ be a $k$--linear small category and let $\C{H}$ be a
  bifunctor on $\C{C}$ with values in $\lmod{k}$.  Define
  $\CH_*(\C{C},\C{H})$ the Hochschild complex of $\C{C}$ with
  coefficients in the bifunctor $\C{H}$ as the differential graded
  $k$--module given by
  \begin{align*}
    \CH_n(\C{C},\C{H})
    := \bigoplus_{X_0,\ldots,X_n} \C{H}(X_0,X_n)\otimes
        {\rm Hom}_\C{C}(X_1,X_0)\otimes\cdots\otimes{\rm Hom}_\C{C}(X_n,X_{n-1})
  \end{align*}
  with a pre-simplicial structure
  \begin{align*}
    \partial_j & (h\otimes X_0\xla{u_1}\cdots\xla{u_n}X_n)\\
     = & \begin{cases}
       \C{H}(u_1,id_{X_n})(h)\otimes X_1\xla{u_2}\cdots\xla{u_n}X_n
            & \text{ if } j=0\\
       h\otimes\cdots\xla{u_{i-1}}X_{i-1}\xla{u_i u_{i+1}}X_{i+1}\xla{u_{i+2}}\cdots
            & \text{ if } 0<j<n-1\\
       \C{H}(id_{X_0},u_n)(h)\otimes X_0\xla{u_1}\cdots\xla{u_{n-1}}X_{n-1}
            & \text{ if } j=n
       \end{cases}
  \end{align*}
  defined for any $n\geq 1$ and $(h\otimes
  X_0\xla{u_1}\cdots\xla{u_n}X_n)$ from $\CH_n(\C{C},\C{H})$.  In the
  case $\C{H}={\rm Hom}_\C{C}$, we denote $\CH_*(\C{C},{\rm
  Hom}_\C{C})$ simply by $\CH_*(\C{C})$.
\end{defn}

Assume $\C{A}_*$ and $\C{B}_*$ are two pre-simplicial $k$--modules
and let $f,g:\C{A}_*\to\C{B}_*$ be two morphisms of pre-simplicial
modules.  Now, recall from \cite{May:SimplicialObjects} that a
pre-simplicial homotopy $h_*$ between $f_*$ and $g_*$ is a set of
$k$--module morphisms $h_i:\C{A}_n\to \C{B}_{n+1}$ defined for
$0\leq i\leq n$ satisfying
\begin{align*}
  h_i\partial_j = & \partial_j h_{i+1},\ \  \text{ if } j\leq i &
  h_i \partial_j = & \partial_{j+1} h_i,\ \ \text{ if } j\geq i+1 &  
  \partial_i h_i = & \partial_i h_{i-1}
\end{align*}
where $f_*=\partial_0 h_0$ and $g_*=\partial_{*+1} h_*$.


\begin{defn}
  Let $\C{C}$ be a $k$--linear category and let $\C{D}$ be a
  $k$--linear subcategory of $\C{C}$.  Then $\C{D}$ is called a
  cofinal subcategory of $\C{C}$ if for every object $C$ of $\C{C}$,
  there exists an object $D$ in $\C{D}$ and a retract $D\xra{r}C$ in
  $\C{C}$.
\end{defn}

\begin{thm}\label{CofinalEquivalence}
  Let $\C{C}$ be a small $k$--linear category and $\C{H}$ be a
  bifunctor on $\C{C}$ with values in $\lmod{k}$.  Assume $\C{D}$ is a
  cofinal subcategory of $\C{C}$.  Then the natural inclusion
  $\CH_*(\C{D},\C{H})\xra{i_*}\CH_*(\C{C},\C{H})$ is a homotopy
  equivalence.
\end{thm}

\begin{proof}
  For every object $C$ in $\C{C}$ fix a choice of object $\delta(C)$
  and a retract $\delta(C)\xra{r(C)}C$ such that for each object $D$
  in $\C{D}$, the choice is $D\xra{id_D}D$.  Denote the left inverse
  of $r(C)$ by $s(C)$ for every $C$ in $\C{C}$.  Now define a morphism
  of differential graded $k$--modules
  $\CH_*(\C{C},\C{H})\xra{M_*}\CH_*(\C{D},\C{H})$ by letting
  \[  M_*(h\otimes X_0\xla{u_1}\cdots\xla{u_n}X_n)
      = \C{H}(s_n,r_0)(h)\otimes \delta(X_0)\xla{s_0 u_1 r_1}\cdots
             \xla{s_{n-1} u_n r_n}\delta(X_n)
  \]
  for any $(h\otimes X_0\xla{u_1}\cdots\xla{u_n}X_n)$ in
  $\CH_n(\C{C},\C{H})$ where we use $r_i = r(X_i)$ and $s_i=s(X_i)$
  for any $0\leq i\leq n$.  It is easy to see that $M_*$ is a morphism
  of pre-simplicial modules since $r_i s_i =id_i$.  Note that the
  composition $M_* i_*$ is identity on $\CH_*(\C{D},\C{H})$.  We are
  going to show $i_* M_*$ is homotopic to the identity on
  $\CH_*(\C{C},\C{H})$.  We define a pre-simplicial homotopy by
  letting
  \begin{align*}
    h_i(h & \otimes X_0\xla{u_1}\cdots\xla{u_n}X_n)\\
     = &\  \C{H}(s_n,id_{X_0})(h)\otimes X_0\xla{u_1}\cdots\xla{u_i}X_i\xla{r_i}
            \delta(X_i)\xla{s_i u_{i+1} r_{i+1}}\cdots\xla{s_{n-1} u_n r_n}\delta(X_n)
  \end{align*}
  for any $0\leq i\leq n$ and for any $(h\otimes
  X_0\xla{u_1}\cdots\xla{u_n}X_n)$ in $\CH_*(\C{C},\C{H})$.  Note that
  $\partial_0 h_0= i_* M_*$ and $\partial_{n+1}h_n = id_*$.  We leave
  the verification of pre-simplicial homotopy identities to the
  reader.
\end{proof}


\begin{cor}
  If $\C{C}'\xra{F}\C{C}$ is a $k$--linear functor, then one has a
  morphism of differential graded $k$--modules of the form
  \[ \CH_*(\C{C}',\C{H}F)\xra{F_*}\CH_*(\C{C},\C{H}) \]  
  for any bifunctor $\C{H}$ on $\C{C}$ with values in $\lmod{k}$.
  Moreover, if $F$ is an equivalence of categories then $F_*$ is a
  homotopy equivalence.
\end{cor}

\begin{proof}
  First, let us explain what $\CH_*(\C{C}',\C{H}F)$ is:
  $\CH_n(\C{C}',\C{H}F)$ generated by homogeneous tensors of the form
  $(h\otimes X_0\xla{u_1}\cdots\xla{u_n}X_n)$ where
  $h\in\C{H}(F(X_n),F(X_0))$ for any $n\geq 0$.  The ``action'' of
  $u_1$ and $u_n$ on $h$ are defined through $F$.  The morphism $F_*$
  of pre-simplicial $k$--modules is defined as
  \[ F_n(h\otimes X_0\xla{u_1}\cdots\xla{u_n}X_n)
     = h\otimes F(X_0)\xla{F(u_1)}\cdots\xla{F(u_n)}F(X_n)
  \]
  for any $n\geq 0$ and $(h\otimes X_0\xla{u_1}\cdots\xla{u_n}X_n)$
  from $\CH_n(\C{C}',\C{H}F)$.  Now assume $F$ is an equivalence with
  a quasi-inverse $G$ and with the isomorphism
  $id_{\C{C}'}\xra{\varphi}GF$.  Note that we have the composition
  \[ \CH_*(\C{C}',\C{H}FGF)\xra{G_*F_*}\CH_*(\C{C}',\C{H}F)
  \]
  and the image of the functor $\C{C}'\xra{GF}\C{C}'$ is a cofinal
  subcategory of $\C{C}'$ since $GF\simeq id_{\C{C}'}$.  Thus
  $G_*F_*\simeq id_*$.  The same argument works also for
  $\C{C}\xra{FG}\C{C}$ and we see that $F_*G_*\simeq id_*$.  The
  result follows.
\end{proof}

\begin{defn}\label{FreeGeneration}
  Let $\C{C}$ be a $k$--linear category which has finite coproducts
  and let $\C{D}$ and $\C{E}$ be two full $k$--linear subcategories.
  $\C{D}$ is said to generate $\C{E}$ if for every object $E$ of
  $\C{E}$ there is a natural number $n\geq 1$ and a set of objects
  $D_1,\ldots,D_n$ of $D$ such that $E\cong \coprod_{i=1}^n D_i$.
\end{defn}

\begin{thm}
  Let $\C{C}$, $\C{D}$ and $\C{E}$ be as in
  Definition~\ref{FreeGeneration}.  Then the natural inclusion
  $\CH_*(\C{D},\C{H})\xra{i_*}\CH_*(\C{E},\C{H})$ is a homotopy
  equivalence for any bifunctor $\C{H}$ on $\C{C}$ with values in
  $\lmod{k}$.
\end{thm}

\begin{proof}
  Take an object $E$ from $\C{E}$ and consider ``$\C{D}$--components''
  $\{D_1,\ldots,D_n\}$ of $E$.  Since $E\cong\coprod_{i=1}^n D_i$ and
  \begin{align*}
    {\rm Hom}_\C{C}(E,E) \cong \bigoplus_{i,j}{\rm Hom}_\C{C}(D_i,D_j)
  \end{align*}
  there are morphisms $E\xla{v_i}D_i$ and $D_j\xla{u_j}E$ such that
  $\sum_iv_iu_i=id_E$.  Now take $h\otimes
  E_0\xla{f_1}\cdots\xla{f_n}E_n$ from $\CH_n(\C{E},\C{H})$ and let
  $D^i_j$ be the $\C{D}$--components of $E_i$, and
  $E_i\xla{v^i_j}D^i_j\xla{u^i_j}E_i$ be the corresponding splitting
  of $id_{E_i}$ for $0\leq i\leq n$.  Define a morphism of
  pre-simplicial modules ${\rm
  CH}_*(\C{E},\C{H})\xra{M_*}\CH_*(\C{D},\C{H})$ by letting
  \begin{align*}
    M_n(h\otimes & E_0\xla{f_1}\cdots\xla{f_n}E_n)\\
      = & \sum_{i_0,\ldots,i_n} \C{H}(u^n_{i_n},v^0_{i_0})(h)
             \otimes D^0_{i_0}\xla{u^0_{i_0}f_1v^1_{i_1}}
             \cdots\xla{u^{n-1}_{i_{n-1}}f_nv^n_{i_n}}D^n_{i_n}
  \end{align*}
  Notice that $M_* i_*$ is the identity on $\CH_*(\C{D},\C{H})$.
  Observe also that the identity $\sum_j v^i_j u^i_j=id_{E_i}$ implies
  $M_*$ is a morphism of pre-simplicial modules.  Although $i_* M_*$
  is not identity, we will furnish a pre-simplicial homotopy between
  $id_*$ and $i_* M_*$ on $\CH_*(\C{E},\C{H})$.  We let
  \begin{align*}
    h_s(h\otimes & E_0\xla{f_1}\cdots\xla{f_n}E_n)\\
    = & \C{H}(id_{E_n},v^0_{i_0})(h)\otimes D^0_{i_0}\xla{u^0_{i_0} f_1 v^1_{i_1}}
            \cdots\xla{u^{s-1}_{i_{s-1}}f_s v^s_{i_s}} D^s_{i_s}
            \xla{u^s_{i_s}}E_s\xla{f_{s+1}}\cdots
            \xla{f_n}E_n
  \end{align*}
  for $n\geq 0$, $0\leq s\leq n$ and $h\otimes
  E_0\xla{f_1}\cdots\xla{f_n}E_n$ from $\CH_n(\C{E},\C{H})$.  We leave
  the verification of the pre-simplicial homotopy identities to the
  reader.
\end{proof}

\section{$B$--categories and equivariant bifunctors}
\label{EquivariantCategories}

\begin{defn}\label{EquivariantCategory}
  A $k$--linear category $\C{C}$ is called a $B$--category if each
  ${\rm Hom}_\C{C}(X,Y)$ is a left $B$--module and the composition
  \[  {\rm Hom}_\C{C}(Y,Z)\times {\rm Hom}_\C{C}(X,Y)\xra{} {\rm Hom}_\C{C}(X,Z) \]
  is a $B$--module morphism via the diagonal action of $B$, for any
  $X,Y,Z$ taken from $Ob(\C{C})$.  In other words, one has
  $b(gf)=b_{(1)}(g)b_{(2)}(f)$ for any $b\in B$, $f\in{\rm
  Hom}_\C{C}(X,Y)$ and $g\in{\rm Hom}_\C{C}(Y,Z)$.  A functor
  $\C{C}\xra{F}\C{C}'$ between two $B$--categories is called
  $B$--equivariant if the structure morphisms
  \[  {\rm Hom}_\C{C}(X,Y)\xra{F_{X,Y}}{\rm Hom}_{\C{C}'}(F(X),F(Y)) \]
  are $B$--module morphisms for any $X,Y\in Ob(\C{C})$.
\end{defn}

\begin{rem}
  It came to our attention that Cibils and Solotar defined the same
  notion in \cite{CibilsSolotar:GoloisCoveringMoritaEquivalence} but
  they called the same object as $B$--module category.  Their primary
  example of the underlying bialgebra is a group ring which is
  cocommutative.  The bialgebras, or in general Hopf algebras, we
  consider are not necessarily cocommutative.
\end{rem}

\begin{exm}\label{BasicExample}
  Assume that $B$ is a Hopf algebra and $A$ is a $B$--module algebra.
  Consider the category $\emod{B}{A}$ of left $B$--equivariant right
  $A$--modules and all $A$--linear morphisms.  Note that we do
  consider {\bf all} $A$--module morphisms, not just $B$--equivariant
  $A$--module morphisms.  Define a left $B$ action of ${\rm
  Hom}_{\emod{B}{A}}(X,Y)$ by letting
  \[  (bf)(x) = b_{(1)}f(S(b_{(2)})x) \]
  for any $f\in {\rm Hom}_{\emod{B}{A}}(X,Y)$ and $x\in X$.
  However, one needs to show that $bf$ is still a right $A$--module
  morphism for any $f\in {\rm Hom}_{\emod{B}{A}}(X,Y)$ and $b\in B$.
  Therefore we check
  \begin{align*}
    (bf)(xa) = &  b_{(1)}f(S(b_{(2)})(xa)) 
             = b_{(1)}f(S(b_{(3)})(x)S(b_{(2)})(a))\\
             = & b_{(1)}f(S(b_{(4)})(x)) b_{(2)}S(b_{(3)})(a)
             = (bf)(x)a
  \end{align*}
  for any $a\in A$, $b\in B$, $x\in X$ and $f\in{\rm
  Hom}_{\emod{B}{A}}(X,Y)$.  Now notice that for any $g\in {\rm
  Hom}_{\emod{B}{A}}(Y,Z)$ one has
  \[  b(gf)(x) = b_{(1)}gf(S(b_{(2)})x)
         = b_{(1)}g(S(b_{(2)})b_{(3)}f(S(b_{(4)})x))
         = (b_{(1)}g)(b_{(2)}f)(x)
  \]
  for any $x\in X$.  In other words ${\emod{B}{A}}$ is a
  $B$--category.
\end{exm}

\begin{exm}
  Assume $B$ is a Hopf algebra and $A$ is a $B$--module algebra.  Let
  $\eproj{B}{A}$ be the full subcategory of $\emod{B}{A}$ consisting
  of finitely generated left $B$--equivariant projective right
  $A$--modules.  Then $\eproj{B}{A}$ is a $B$--category.
\end{exm}

\begin{exm}
  Assume $B$ is a Hopf algebra and $A$ is a $B$--module algebra.
  Consider the full subcategory $*_B^A$ of $\efree{B}{A}$ which
  consists of one single object $A$ considered as a right $A$--module
  via the right regular representation.  Then $*_B^A$ is a
  $B$--category.
\end{exm}

\begin{defn}
  A bifunctor $\C{H}$ on a $B$--category $\C{C}$ with values in
  $\lmod{B}$ is called a $B$--equivariant bifunctor if the structure
  morphisms
  \begin{align}\label{BifunctorStructure}
    \C{H}(X,Y)\otimes & {\rm Hom}_\C{C}(Y,Z)\xra{} \C{H}(X,Z) & 
    {\rm Hom}_\C{C}(W,X)\otimes & \C{H}(X,Y)\xra{} \C{H}(W,Y)
  \end{align}
  are $B$--module morphisms where $B$ acts diagonally on the left.  In
  other words 
  \[ b(\C{H}(u,v)(h)) = \C{H}(b_{(1)}(u), b_{(3)}(v))(b_{(2)}(h)) \]
  for any $u\in Hom_\C{C}(W,X)$, $v\in {\rm Hom}_\C{C}(Y,Z)$ and $h\in
  \C{H}(X,Y)$, and $b\in B$.
\end{defn}

\begin{exm}
  For a $B$--category $\C{C}$, the bifunctor $\C{H}={\rm
  Hom}_\C{C}(\cdot,\cdot)$ is a $B$--equivariant bifunctor on $\C{C}$
  with values in $\lmod{B}$.
\end{exm}

\begin{defn}
  For a $B$--category $\C{C}$, let ${}^B\C{C}$ denote the subcategory
  of morphisms of $\C{C}$ which are $B$--invariant, i.e. $X\xra{f}Y$
  belongs to ${}^B\C{C}$ iff $b(f)=\varepsilon(b)f$ for any $b\in B$.
\end{defn}

\begin{rem}\label{NaiveQuotient}
  Assume $B$ is a Hopf algebra and ${\emod{B}{A}}$ be the
  $B$--category defined in Example~\ref{BasicExample}.  Then one can
  see that ${}^B\emod{B}{A}$, the subcategory of $B$--invariant
  morphisms, is the category of left $B$--equivariant right
  $A$--modules and their $B$--equivariant $A$--module morphisms since
  $bf=\varepsilon(b)f$ iff $f$ is $B$--equivariant.  We would like to
  note that for the $B$--equivariant homotopical invariants of
  ${\emod{B}{A}}$ the subcategory ${}^B\emod{B}{A}$ is rather small.
  The situation is very similar to topological spaces admitting an
  action of a fixed group $G$.  The $G$--equivariant homotopical
  invariants of a $G$--space $X$ are computed via
  $\C{B}(G,X):=EG\underset{G}{\wedge} X$ rather than $X/G\simeq
  *\underset{G}{\wedge} X$.  Similarly, ${}^B\emod{B}{A}$ should be
  considered as the lowest order equivariant invariant of
  $\emod{B}{A}$.  Thus we propose that for higher order equivariant
  homotopical invariants of a $B$--module algebra $A$, such as
  equivariant $K$--theoretical, Hochschild and cyclic homological
  invariants, one should use the $B$--category $\emod{B}{A}$ of
  $B$--equivariant modules and their $A$--linear morphisms, or its
  various $B$--subcategories, instead of using simply
  ${}^B\emod{B}{A}$ the subcategory of $B$--equivariant $A$--module
  morphisms of $\emod{B}{A}$.  Our justification lies in
  Section~\ref{Morita} where we prove Morita invariance in
  Corollary~\ref{EquivariantMorita}.
\end{rem}

\section{Morita invariance}
\label{Morita}

\begin{lem}
  If $\C{C}$ is a $B$--category and $\C{H}$ is a bifunctor on $\C{C}$
  with values in $\lmod{B}$, then $\CH_*(\C{C},\C{H})$ is a graded
  $B$--module.  However, $\CH_*(\C{C},\C{H})$ is not a pre-simplicial
  $B$--module unless $B$ is cocommutative.
\end{lem}

\begin{proof}
  The $B$--action of $\CH_*(\C{C},\C{H})$ is defined diagonally, i.e.
  \[ L_b(h\otimes X_0\xla{u_1}\cdots\xla{u_n}X_n)
     := b_{(1)}(h)\otimes X_0\xla{b_{(2)}u_1}\cdots\xla{b_{(n+1)}u_n}X_n
  \] 
  for any $b\in B$, $n\geq 0$ and $h\otimes
  X_0\xla{u_1}\cdots\xla{u_n}X_n$ from $\CH_n(\C{C},\C{H})$.  The fact
  that $\CH_*(\C{C},\C{H})$ is NOT a pre-simplicial $B$--module is
  because of the last face morphism: one has
  \begin{align*}
     \partial_n L_b & (h\otimes X_0\xla{u_1}\cdots\xla{u_n}X_n)\\
     = & \C{H}(b_{(n+1)}(f_n),id_{X_0})(b_{(1)}h)\otimes 
        X_0\xla{b_{(2)}u_1}\cdots\xla{b_{(n)}u_{n-1}}X_{n-1}
  \end{align*}
  which is different than
  \begin{align*}
    L_b\partial_n & (h\otimes X_0\xla{u_1}\cdots\xla{u_n}X_n)\\
     = & \C{H}(b_{(1)}(f_n),id_{X_0})(b_{(2)}h)\otimes 
        X_0\xla{b_{(3)}u_1}\cdots\xla{b_{(n+1)}u_{n-1}}X_{n-1}
  \end{align*}
  for any $n\geq 0$, $b\in B$ and $h\otimes
  X_0\xla{u_1}\cdots\xla{u_n}X_n$ from $\CH_*(\C{C},\C{H})$, unless
  $B$ is cocommutative.
\end{proof}

\begin{defn}
  Define a graded $k$--submodule $J_*(\C{C},B,\C{H})$ of
  $\CH_*(\C{C},\C{H})$ generated by elements of the form
  \[ [L_b,\partial_n](h\otimes X_0\xla{u_1}\cdots\xla{u_n}X_n)  \]
  where $n\geq 1$, $b\in B$ and $h\otimes
  X_0\xla{u_1}\cdots\xla{u_n}X_n$ from $\CH_n(\C{C},\C{H})$.
\end{defn}

\begin{lem}
  $J_*(\C{C},B,\C{H})$ is a differential graded $k$--submodule and
  graded $B$--submodule of $\CH_*(\C{C},\C{H})$.  Therefore
  $\CH_*(\C{C},\C{H})/J_*(\C{C},B,\C{H})$ is a differential graded
  $B$--module.
\end{lem}

\begin{proof}
  The proof is identical to that of Proposition~\ref{Quotient}.
\end{proof}

\begin{defn}
  Let $\QCH_*(\C{C},B,\C{H})$ be the quotient
  $\CH_*(\C{C},\C{H})/J_*(\C{C},B,\C{H})$.
\end{defn}

\begin{thm}
  Let $\C{D}$ be a cofinal $B$--subcategory of $\C{C}$.  Then for any
  $B$--equivariant bifunctor $\C{H}$ on $\C{C}$ with values in
  $\lmod{B}$ then the natural inclusion
  \[ \QCH_*(\C{D},B,\C{H})\xra{i_*}\QCH_*(\C{C},B,\C{H}) \]
  is a homotopy equivalence.
\end{thm}

\begin{proof}
  The proof is almost verbatim that of
  Theorem~\ref{CofinalEquivalence} after noticing $J_*(\C{C},B,\C{H})$
  is stable under the pre-simplicial homotopy we furnished there.
\end{proof}

\begin{cor}
  Let $\C{C}$ and $\C{C}'$ be two $B$--categories and let
  $\C{C}\xra{F}\C{C}'$ be a functor of $B$--categories.  Then for any
  $B$--equivariant bifunctor $\C{H}$ on $\C{C}'$ one has a morphism of
  differential graded $B$--modules of the form
  \[  \QCH_*(\C{C},B,\C{H}F)\xra{F_*}\QCH_*(\C{C}',B,\C{H})  \]
  Moreover, if $F$ is an equivalence of $B$--categories, then $F_*$ is
  a homotopy equivalence.
\end{cor}

\begin{cor}
  Assume $A$ is a $B$--module algebra.  Then $\efree{B}{A}$ is a
  cofinal $B$--subcategory of $\eproj{B}{A}$.  Furthermore, the
  natural inclusion functor $\efree{B}{A}\xra{}\eproj{B}{A}$ induces a
  homotopy equivalence of differential graded $B$--modules of the form
  \[ \QCH_*(\efree{B}{A},B,\C{H}) \xra{i_*}\QCH_*(\eproj{B}{A},B,\C{H}) \]
  for any $B$--equivariant bifunctor $\C{H}$ on $\eproj{B}{A}$ with
  values in $\lmod{B}$.
\end{cor}

\begin{cor}
  Assume $A$ is a $B$--module algebra.  The subcategory $*_B^A$ freely
  generates the $B$--subcategory $\efree{B}{A}$ of $\eproj{B}{A}$.
  Then the natural inclusion
  \begin{align*}
    \QCH_*(*_B^A,B,\C{H})\xra{i_*}\QCH_*(\efree{B}{A},B,\C{H})
  \end{align*}
  is a homotopy equivalence of differential graded $B$--modules for
  any $B$--equivariant bifunctor $\C{H}$ defined on $\eproj{B}{A}$
  with values in $\lmod{B}$.  Furthermore, $\QCH_*(*_B^A,B,\C{H})$ and
  $\QCH_*(A,B,\C{H}(A,A))$ are isomorphic as the differential graded
  $B$--modules.
\end{cor}

\begin{proof}
  The proof relies on the fact that ${\rm Hom}_{\eproj{B}{A}}(A,A)\cong
  A$.  The rest of the proof is trivial.
\end{proof}

\begin{thm}[Morita invariance for Hopf--Hochschild (co)homology]\label{EquivariantMorita}
  Let $B$ be a Hopf algebra and assume that $A$ and $A'$ are two
  $B$--module algebras.  If $\emod{B}{A}$ and $\emod{B}{A'}$ the
  category of finitely generated $B$--equivariant representations of
  $A$ and $A'$ respectively are $B$--equivariantly equivalent, then
  Hopf--Hochschild complex ${}_B\QCH_*(A,B,A)$ of $A$ and
  Hopf--Hochschild complex ${}_B\QCH_*(A',B,A')$ of $A'$ are
  quasi-isomorphic.
\end{thm}

\section{Twisted equivariant bifunctors as coefficients}
\label{TwistedBifunctors}

In this section we assume $B$ is a Hopf algebra with an invertible
antipode.

\begin{defn}
  An arbitrary left-left $B$--module/comodule $M$ is called a Yetter
  Drinfeld module if one has
  \begin{align*}
    (bm)_{(-1)}\otimes (bm)_{(0)} = b_{(1)}m_{(-1)}S(b_{(3)})\otimes b_{(2)}m_{(0)}    
  \end{align*}
  for any $b\in B$ and $m\in M$ where we use Sweedler's notation for
  the coproduct on $H$ and for the $H$--coaction on $M$.
\end{defn}

\begin{defn}\label{Twisting}
  Assume $M$ is a Yetter-Drinfeld module over $B$.  Let $\C{C}$ be an
  arbitrary $B$--category and $\C{H}$ be an arbitrary $B$--equivariant
  bifunctor on $\C{C}$ with values in $\lmod{B}$.  We define a new
  bifunctor $M\ltimes \C{H}$ by letting $M\ltimes\C{H}(X,Y):=M\otimes
  \C{H}(X,Y)$ on the objects for any $X,Y\in Ob(\C{C})$.  Now we let
  \[ b(m\otimes h) = b_{(1)}(m)\otimes b_{(2)}(h) \]
  for any $m\otimes h$ in $\C{H}(X,Y)$.  We define the bifunctor
  $M\ltimes\C{H}$ on morphisms as follows: notice that such functors are
  defined by structure morphisms given in~(\ref{BifunctorStructure}).
  Then we let
  \[  (m\otimes h)(Y\xla{\beta} Z) := m\otimes \C{H}(id_X,\beta)(h) \]
  and
  \[  (W\xla{\alpha}X)(m\otimes h) := m_{(0)}\otimes\C{H}(S^{-1}(m_{(-1)})(\alpha),id_Y)(h) \]
  for any $m\otimes h$ from $M\ltimes\C{H}(X,Y)$, $\alpha\in {\rm
  Hom}_\C{C}(W,X)$ and $\beta\in{\rm Hom}_\C{C}(Y,Z)$.  For simplicity
  we will denote $\C{H}(\alpha,\beta)(h)$ by $\alpha h \beta$ for any
  $h\in\C{H}(X,Y)$, $\alpha\in{\rm Hom}_\C{C}(W,X)$ and $\beta\in{\rm
  Hom}_\C{C}(Y,Z)$.
\end{defn}

\begin{lem}
  Let $B$, $M$, $\C{C}$ and $\C{H}$ be as in
  Definition~\ref{Twisting}.  Then $M\ltimes\C{H}$ is also a
  $B$--equivariant bifunctor on $\C{C}$ with values in $\lmod{B}$.
\end{lem}

\begin{proof}
  For $b\in B$ and $m\otimes h$ in $M\ltimes\C{H}(X,Y)$, and $\alpha\in{\rm
  Hom}_\C{C}(W,X)$ and ${\rm Hom}_\C{C}(Y,Z)$ we consider
  \[ b((m\otimes h)\beta) = b(m\otimes h\beta) 
     = b_{(1)}(m)\otimes b_{(2)}(h\beta)
     = b_{(1)}(m)\otimes b_{(2)}(h)(b_{(3)}\beta)
     = b_{(1)}(m\otimes h)(b_{(2)}\beta)
  \]
  and
  \begin{align*}
    b(\alpha(m\otimes h))
     = & b(m_{(0)}\otimes (S^{-1}(m_{(-1)})\alpha)(h))\\
     = & b_{(1)}(m_{(0)})\otimes (b_{(2)}S^{-1}(m_{(-1)})\alpha)b_{(3)}(h)\\
     = & b_{(2)(2)}m_{(0)}\otimes (b_{(2)(3)}S^{-1}(m_{(-1)})S^{-1}(b_{(2)(1)}) b_{(1)}\alpha)
               b_{(3)}(h)\\
     = & (b_{(1)}\alpha) (b_{(2)}m\otimes b_{(3)}h)
     = (b_{(1)}\alpha) b_{(2)}(m\otimes h)
  \end{align*}
  which implies $M\ltimes\C{H}$ is a $B$--equivariant bifunctor on
  $\C{C}$ with values in $\lmod{B}$ whenever $\C{H}$ is a bifunctor on
  $\C{C}$ with values in $\lmod{B}$ and $M$ is a Yetter-Drinfeld
  module on $B$.
\end{proof}

For a $B$--module algebra $A$ and a Yetter--Drinfeld module $M$, we
now define a new differential graded $B$--module
$\QCH_*(A,B,A;M):=\QCH_*(*^A_B,B,M\ltimes {\rm Hom}_A(\cdot,\cdot))$
and obtain twisted version of the Morita equivalence

\begin{thm}[Morita invariance of twisted Hopf--Hochschild (co)homology]\label{TwistedMorita}
  Let $B$ be a Hopf algebra and $M$ be a Yetter-Drinfeld module over
  $B$.  Assume also that $A$ and $A'$ are two $B$--module algebras.
  If the category of finitely generated $B$--equivariant
  representations of $A$ and $A'$ are $B$--equivariantly equivalent,
  then the twisted Hopf--Hochschild complexes ${}_B\QCH_*(A,B,A;M)$
  and ${}_B\QCH_*(A',B,A';M)$ are quasi-isomorphic.
\end{thm}

\bibliographystyle{plain} 
\bibliography{bibliography}

\begin{thebibliography}{10}

\bibitem{KhalkhaliAkbarpour:EquivariantCyclicHomologyII}
R.~Akbarpour and M.~Khalkhali.
\newblock Equivariant cyclic cohomology of {$H$}-algebras.
\newblock {\em $K$--Theory}, 29(4):231--252, 2003.

\bibitem{AtiyahSegal:EquivariantKTheory}
M.~F. Atiyah and G.~B. Segal.
\newblock Equivariant {$K$}-theory and completion.
\newblock {\em J. Differential Geometry}, 3:1--18, 1969.

\bibitem{BaajSkandalis:HopfEquivariantKKTheory}
S.~Baaj and G.~Skandalis.
\newblock {$C^*$}-alg\`ebres de {H}opf et th\'eorie de {K}asparov
  \'equivariante.
\newblock {\em $K$-Theory}, 2(6):683--721, 1989.

\bibitem{BlockGetzlerJones:CyclicHomologyOfCrossedProductAlgebras}
J.~Block, E.~Getzler, and J.~D.~S. Jones.
\newblock The cyclic homology of crossed product algebras. {II}. {T}opological
  algebras.
\newblock {\em J. Reine Angew. Math.}, 466:19--25, 1995.

\bibitem{Brylinski:CyclicHomologyAndEquivariantTheories}
J.-L. Brylinski.
\newblock Cyclic homology and equivariant theories.
\newblock {\em Ann. Inst. Fourier (Grenoble)}, 37(4):15--28, 1987.

\bibitem{CibilsSolotar:GoloisCoveringMoritaEquivalence}
C.~Cibils and A.~Solotar.
\newblock Galois coverings, {M}orita equivalence and smash extensions of
  categories over a field.
\newblock {\em Documenta Math.}, (11):143--159, 2005.

\bibitem{ConnesMoscovici:HopfCyclicCohomology}
A.~Connes and H.~Moscovici.
\newblock Hopf algebras, cyclic cohomology and transverse index theorem.
\newblock {\em Comm. Math. Phys.}, 198:199--246, 1998.

\bibitem{GetzlerJones:CyclicHomologyOfCrossedProductAlgebras}
E.~Getzler and J.~D.~S. Jones.
\newblock The cyclic homology of crossed product algebras.
\newblock {\em J. Reine Angew. Math.}, 445:161--174, 1993.

\bibitem{Green:KTheoryOfCrossedProductAlgebras}
P.~Green.
\newblock Equivariant {$K$}-theory and crosed product {$C^*$}-algebras.
\newblock Number~38, pages 337--338, Providence, 1982. Amer. Math. Soc.

\bibitem{Khalkhali:HopfCyclicHomology}
P.~M. Hajac, M.~Khalkhali, B.~Rangipour, and Y.~Sommerh\"auser.
\newblock Hopf--cyclic homology and cohomology with coefficients.
\newblock {\em C. R. Math. Acad. Sci. Paris}, 338(9):667--672, 2004.

\bibitem{Julg:KTheoryOfCrossedProductAlgebras}
P.~Julg.
\newblock {$K$}-th\'eorie \'equivarante et produits crois\'es.
\newblock {\em C.R. Acad. Sci. Paris}, 292(13):629--632, 1981.

\bibitem{Kasparov:EquivariantKKTheory}
G.~G. Kasparov.
\newblock Equivariant {$KK$}-theory and the {N}ovikov conjecture.
\newblock {\em Invent. Math.}, 91(1):147--201, 1988.

\bibitem{Kaygun:BialgebraCyclicK}
A.~Kaygun.
\newblock {Bialgebra cyclic homology with coefficients}.
\newblock {\em $K$--Theory}, 34(2):151--194, 2005.

\bibitem{Kaygun:BivariantHopf}
A.~Kaygun and M.~Khalkhali.
\newblock Bivariant {H}opf cyclic cohomology.
\newblock Preprint at {\tt arXiv:math.KT/0606341}.

\bibitem{May:SimplicialObjects}
J.~P. May.
\newblock {\em Simplicial objects in algebraic topology}.
\newblock University of Chicago Press, 1966.

\bibitem{McCarthy:CyclicHomologyOfExactCategories}
R.~McCarthy.
\newblock The cyclic homology of an exact category.
\newblock {\em J. Pure Appl. Algebra}, 93(3):251--296, 1994.

\bibitem{NeshveyevTuset:HopfEquivariantKTheory}
S.~Neshveyev and L.~Tuset.
\newblock Equivariant cyclic cohomology, {$K$}--theory and index formulas.
\newblock {\em $K$--Theory}, (31):357--378, 2004.

\bibitem{Pirashvili:HochschildAndCyclicHomologyAsFunctorHomology}
T.~Pirashvili and B.~Richter.
\newblock Hochschild and cyclic homology via functor homology.
\newblock {\em $K$-Theory}, 25(1):39--49, 2002.

\end{thebibliography}

\end{document}